\documentstyle{amsart}

\newtheorem{theo}{Theorem}

\newtheorem{coro}{Corollary}
\newtheorem{exam}{Example}

\begin{document}
\title[Classes of forms Witt equivalent to a second trace form of fields]
{Classes of forms Witt equivalent to a second trace form over
fields of characteristic two}
\author{Ana Cecilia de la Maza}
\address{Instituto de Matem\'atica y F\'{\i}sica,
 Universidad de Talca, Casilla 747, Talca, Chile}
\thanks{This work  was supported by Fondecyt grants N$^{o} 1010205$, and by
Programa formas extremas y representaci\'on de formas
cuadr\'aticas, Universidad de Talca.}

\subjclass{ Primary 11 E 04; Secondary 11 E 81}


\begin{abstract}
Let $F$ be a field of characteristic two. We determine all
non-hyperbolic quadratic forms over $F$ that are Witt equivalent
to a second trace form.
\end{abstract}

\maketitle
\section{Introduction}
 Let $E/F$ be a finite separable field extension.
We define the trace form for this extension by
$q(x)=tr_{E/F}(x^2)$. When the characteristic of $F$ is not equal
to $2$, the trace form $(E,q)$ is non-degenerate. However, if the
characteristic is $2$ then $(E,q)$ is degenerate and splits as
$[1]\perp V$, with $V$ totally isotropic. It is therefore natural
to introduce a modified ``second trace form".  To this end one
considers for each $a\in E$ its characteristic polynomial
\begin{equation}\label{polc}
p(x,a)=x^n-T_1(a)x^{n-1}+T_2(a)x^{n-2}+\cdots+(-1)^{n}T_n(a)
\end{equation}
(whence $T_1(a)=tr_{E/F}(a)$ and $T_n(a)$ is the norm of $a$). It
is clear that $(E,T_2)$ is a quadratic form. When the degree $n$
of the extension is odd this form is necessarily singular. To
arrive at a non-degenerate form, two methods have been proposed in
the literature. One method, due to Berg\'e and Martinet \cite{bm},
 increases the dimension of the space by $1$ using the \'etale
$F$-algebra. The other method, due to Revoy \cite{r},  reduces the
dimension of the space by $1$. In this note we will adopt the
second method and call such forms  {\it 2-trace forms}.

We consider the following problem: Which elements $[q]\neq 0$ of
the Witt-group $W_q(F)$ are represented by 2-trace forms?. Our
theorem 3 fully answers this question. Moreover, we will partially
answer the same question for $[q]=0$ (see Prop. 1 and Prop. 2).
For fields of characteristic not equal to $2$ this problem seems quite
more complicated: for partial results concerning generic fields one may consult
\cite{cp} and \cite{ehp}; a complete solution for Hilbertian fields is given in
 \cite{sch}, \cite{ks} and \cite{wat}.

\section{The second trace form}

As we remarked in the introduction, there are two ways to define a
second trace form. In this section we will prove that in fact the
corresponding forms are Witt equivalent.

Let $E/F$ be a finite separable field extension. The second trace
form $T_{E/F}$ of the extension $E/F$ was defined by Revoy
\cite{r} as $(E,T_2)$ if the degree $[E:F]$  is even, and as
$(E_0,T_2)$ if the degree is odd, where $T_1,\ T_2$ are given by
\eqref{polc} and $E_0=Ker\ T_1$.  It is important to remark that
the bilinear form $b_q$ associated to $T_{E/F}$ satisfies the
following relations:
\begin{eqnarray}
b_q(x,y)=T_2(x+y)-T_2(x)-T_2(y)=T_1(xy)-T_1(x)T_1(y), \label{bq}\\
 T_1(x^2)=(T_1(x))^2\quad \text{ and } \quad
b_q(x^2,y^2)=b_q(x,y)^2.\label{t1t2}
\end{eqnarray}

On the other hand, Berg\'e and Martinet defined in \cite{bm}  the
second trace form as Revoy if the degree $[E:F]$ is even and as
$(E\times F,T_2)$ if not.

\begin{theo} The Revoy form and the Berg\'e-Martinet form are Witt
equivalent.
\end{theo}
\noindent Proof. If $[E:F]$ is odd, then the form $(E\times
F,T_2)$ of Berg\'e and Martinet \cite[p. 14]{bm} splits  as
follows $(F (1,0)+ F(0,1))\perp E_0\times \{0\}$. Since $F (1,0)+
F(0,1)$ is an hyperbolic plane, the claim follows immediately.
\hfill $\Box$
\section{2-algebraic forms}
In this section we determine all non-hyperbolic quadratic forms
over $F$ that are Witt equivalent to some second trace form.
Furthermore we give fields where hyperbolic forms are Witt
equivalent to a second trace form.

\begin{theo}
Let $E/F$ be a finite  separable field extension, with
$[E:F]=2n+1$ or $[E:F]=2n$. Then $T_{E/F}=(n-1){\mathbb{H}}\perp
[1,a]$, for some $a\in F$.
\end{theo}

\indent Proof. The assertion is deduced from Theorem 1 above and
Theorem 3.5 in \cite[p. 13-14]{bm}. In order to illustrate the
ideas, we will give the proof in case $[E:F]=2n+1$.

Since $(E_0,T_2)$ is nonsingular, there exists a symplectic basis
$\{e_i,f_i\}_{1\leq i\leq n}$ of $E_0$. Hence for each $1\leq
i\leq n$ there exists $x_i, \ y_i\in F e_i+F f_i$ such that
$T_2(x_i)\neq 0$ and $b_q(x_i,y_i)=1$ (note that $[0,0]=[1,0]$
\cite[p. 150]{sa1}).

Put $e'_i:=T_2(x_i)^{-1}x_i^2$ and $f'_i:=T_2(x_i)y_i^2$. We have
by \eqref{t1t2} that $e'_i,\ f'_i\in E_0$ and $b_q(e'_i,f'_i)=1$.
Since $T_2(e'_i)=1$,  we obtain a symplectic basis such that the
quadratic form decomposes as follows
\[T_{E/F}= [1,a_1]+[1,a_2]+\cdots + [1,a_n],\]
where $a_i=T_2(f'_i)$. Using the relation
$[1,b]+[1,c]\cong[1,b+c]+\mathbb{H}$ over $F$, we obtain
$T_{E/F}\cong (n-1){\mathbb{H}}+[1,a]$ with $a\equiv
a_1+a_2+\cdots a_n \mod\ \wp(F)$, where $\wp(F):=\{x^2+x\mid x\in
F\}$. \hfill $\Box$ \vspace{0.1in}

We call {\it 2-algebraic} a form over $F$ if it is  Witt
equivalent to some second trace form.

\begin{coro}
A   non  hyperbolic quadratic form $(V,q)$ over $F$ is 2-algebraic
 if and only if $(V,q)=r{\mathbb{H}}\perp (V_a,q_a),$
with $V_a$ an anisotropic plane representing  1.
\end{coro}

\noindent Proof: $\Rightarrow$) Is clear by  Theorem 2.

\noindent $\Leftarrow$) Let $(V,q) =r{\mathbb{H}} \perp
(V_a,q_a)$, with $(V_a,q_a)$ 2-dimensional anisotropic
representing  1. Then we can rewrite $(V_a,q_a)=[1,b]$, where
$b\notin \wp(F)$. Let $E=F(\alpha)$ with $\alpha \in \overline{F}$
and $\alpha^2+\alpha+b=0$. We have that $(E,T_{E/F})=[1,b]$.
\hfill $\Box$

\begin{exam} Let $F={\mathbb{F}}_2(a)$ and $E=F(b)$, where
 $a^2+a+1=0$ and $b^3+b+a=0$. Then $T_{E/F}=[1,a]\neq [1,1]$.

In fact, using \eqref{bq}, we see that $\{b^2,(1+a)b\}$ is a
symplectic basis for $E_0= Ker\ T_1$. Since
$p(x,(1+a)b)=x^3+ax+a$, $1\in \wp(F)$  and  $a\notin \wp(F)$, we
obtain  the form $(E_0,T_2)=[1,a]\neq [1,1]$.
\end{exam}

\begin{coro}
If a non  hyperbolic quadratic form $(V,q)$ over $F$ is
2-algebraic then there exists a quadratic extension field $E$ of $F$
such that the extension $(V\otimes_FE,q_E)$ is hyperbolic.
\end{coro}

\noindent Proof: See the proof of Theorem3 and note that
$[1,b]=\mathbb{H}$ over $E=F(\alpha)$, with $\alpha^2+\alpha+b=0$
(see \cite[p. 150]{sa1})  . \hfill $\Box$

\begin{theo} Let $F={\mathbb{F}}_2$  or  $F={ \mathbb{F}}_2(t)$
 with $t$ transcendental over ${ \mathbb{F}}_2$. Then
hyperbolic quadratics form over $F$ are 2-algebraic.
\end{theo}

\indent Proof. We only need to find an extension $E$ of $F$ such
that $T_{E/F}$ is hyperbolic. We first remark that
$p(x):=x^4+x^3+1$ is irreducible over ${\mathbb{F}}_2$ and also
over ${\mathbb{F}}_2(t)$. Let $\alpha$ be a root of $p$ and
$E=F(\alpha)$. We decompose the trace form $T_{E/F}$ with
respect to the basis $\{\alpha,\ 1+\alpha^3\}\cup \{\alpha^2,
\alpha +\alpha^2+\alpha^3\}$. Noting that this basis has the
elements  conjugate to $\alpha$, it is easy to recognise that each
vector basis is isotropic, and furthermore by \eqref{bq} we see
that it is a symplectic basis. Hence, the space is hyperbolic.
\hfill $\Box$

\begin{theo} Let $F$ be  a field.
If there exists $a\in F^*$ and $n$ odd such that the polynomial
$x^n-a$ is irreducible over $F[x]$, then  hyperbolic quadratics
space over $F$ are 2-algebraic.
\end{theo}

\noindent Proof. Let $E=F(\alpha)$, where $\alpha\in \overline{F}$
and   $\alpha^n=a$. For $1\leq k\leq n-1$, the linear
transformation $f_{\alpha^k}: x\mapsto x\alpha^k$ is given by the
matrix $c_{ij}(k)$, where
\[c_{ij}(k)=\left\{\begin{array}{ll}
1&\mbox{ if }  j=i-k\\
a&\mbox{ if }  j=n+i-k\\
0& \mbox{ otherwise }\end{array}\right.
\]
Then for $1\leq k\leq n-1$,  $\alpha^k\in E_0$, because
$c_{ii}(k)=0$ for each $i$. Noting that $T_1(a)=a$ we obtain the
decomposition
\[E_0=\langle \alpha,\alpha^{n-1}\rangle \perp \langle \alpha^2,\alpha^{n-2}\rangle
\perp \cdots \perp \langle \alpha^{
\frac{n-1}{2}},\alpha^{\frac{n+1}{2}}\rangle,\] where $\langle
x,y\rangle$ is the  space generate by $x$ and $y$.
 Hence, using that  $n\neq 2k$, we deduce  that $T_2(\alpha^k)=0$
for  $1\leq k\leq n-1$,so $(E_0,T_2)=(\frac{n-1}{2})\mathbb{H}.$
\hfill $\Box$

\end{document}